\documentclass[10 pt]{article}

\usepackage{latexsym}
\usepackage{amssymb}

\title{Markov Extensions for Dynamical Systems with Holes: An Application to Expanding Maps of the Interval}
\author{Mark F. Demers\thanks{Department of Mathematics, Georgia Institute of Technology, Atlanta GA 30332.
demers@math.gatech.edu.}}

\begin{document}

\newtheorem{theorem}{Theorem}[section]
\newtheorem{definition}[theorem]{Definition}
\newtheorem{proposition}[theorem]{Proposition}
\newtheorem{claim}[theorem]{Claim}
\newtheorem{lemma}[theorem]{Lemma}

\newcommand{\T}{\hat{T}}
\newcommand{\M}{\hat{M}}
\newcommand{\dklj}{\Delta_{l,j}^{(k)}}
\newcommand{\dkljs}{\Delta_{l,j}^{(k)*}}
\newcommand{\dilj}{\Delta_{l,j}^{(i)}}
\newcommand{\diljs}{\Delta_{l,j}^{(i)*}}
\newcommand{\dlj}{\ensuremath{\Delta_{l,j}}}
\newcommand{\vone}{\ensuremath{\varphi_1}}
\newcommand{\vtwo}{\ensuremath{\varphi_2}}
\newcommand{\lone}{\ensuremath{\lambda_1}}
\newcommand{\ltwo}{\ensuremath{\lambda_2}}
\newcommand{\di}{\ensuremath{\Delta^{(i)}}}
\newcommand{\hdi}{\ensuremath{\hat{\Delta}^{(i)}}}
\newcommand{\invx}{\ensuremath{\displaystyle 
    \sum_{\stackrel{y \in \dkljs}{\scriptscriptstyle F(y)=x}}}}
\newcommand{\invz}{\ensuremath{\displaystyle 
    \sum_{\stackrel{w \in \dkljs}{\scriptscriptstyle F(w)=z}}}}
\newcommand{\invxz}{\ensuremath{\displaystyle 
    \sum_{\stackrel{y,w \in \dkljs}{\scriptscriptstyle F(y)=x, F(w)=z}}}}
\newcommand{\pf}{\ensuremath{\mathcal{P}f}}
\newcommand{\hilj}{\tilde{H}_{l,j}^{(i)}}
\newcommand{\Z}{\ensuremath{\mathcal{Z}}}
\newcommand{\done}{\ensuremath{\Delta^{(1)}}}
\newcommand{\C}{\ensuremath{\tilde{C}}}
\newcommand{\alf}{\ensuremath{^{\alpha}}}
\newcommand{\h}{\ensuremath{\tilde{H}}}
\newcommand{\Q}{\ensuremath{\mathcal{Q}}}
\newcommand{\Li}{\ensuremath{\Lambda^{(i)}}}
\newcommand{\F}{\hat{F}}
\newcommand{\pilj}{\ensuremath{\pi_{l,j}^{(i)-1}}}
\newcommand{\I}{\hat{I}}

\maketitle

\begin{abstract}
We introduce the Markov extension, represented schematically as a tower, to the study of dynamical systems with holes. 
For tower maps with small holes, we prove the existence of 
conditionally invariant probability measures which are absolutely
continuous with respect to Lebesgue measure (abbreviated a.c.c.i.m.).
We develop restrictions on the Lebesgue measure of the holes and simple
conditions on the dynamics of the tower which ensure existence
and uniqueness in a class of Holder continuous densities.  
We then use these results
to study the existence and properties of a.c.c.i.m. for $C^{1+ \alpha}$ expanding maps of the interval with holes.
We obtain the convergence of the a.c.c.i.m. to the SRB measure of the corresponding closed system
as the measure of the hole shrinks to zero.
\end{abstract}

\section{Introduction}
\label{introduction}

The study of dynamical systems with holes was launched by Pianigiani and Yorke in their
seminal paper \cite{pianigiani yorke}.  
In it, they posed the following still open question.

Consider a particle on a billiard table with convex boundaries so that the dynamics of the particle are
hyperbolic, i.e.\ the trajectories are unstable with respect to initial conditions.  
Suppose a small hole is made in the table.  
What are the statistical properties of the trajectories in this system?  
If $p_n$ is the probability that a trajectory remains on the table until time $n$, what is the decay
rate of $p_n$?  
More generally, we can place a particle randomly on the table according to an initial distribution $\mu_0$.  
If $\mu_n$ represents its normalized distribution at time $n$ (assuming 
the particle has not escaped by time $n$), does $\mu_n$ converge to some $\mu$ independent of $\mu_0$?
Such a measure $\mu$ is a {\em conditionally invariant measure} for the open billiard system.

Considering the billiard table with a small hole as a perturbation 
of the billiard table with no holes, we can pose a related question in terms of the stability of the
closed system:  does the conditionally invariant measure of the open system converge to the invariant
measure of the closed system as the size of the hole tends to zero?

Although these questions remain open, much progress has been made in understanding the existence
and properties of conditionally invariant measures for dynamical systems with holes.  
More generally, the problem can be stated as follows.

Let $\T$ be a piecewise differentiable map of a Riemannian manifold $\M$.  
We take the hole to be an open set $H$ in $\M$ and keep track of the iterates of points until they reach the hole.  
Once a point enters $H$, it is not allowed to return.

Let $M= \M \backslash H$ and let $T = \T|(M \cap \T^{-1}M)$.   
A probability measure $\mu$ on $M$ is said to be a conditionally invariant measure with respect to $T$ if $\mu$ satisfies
\begin{equation}
\label{conditional invariance}
\frac{\mu (T^{-1}A)}{\mu (T^{-1}M)} = \mu(A) 
\end{equation}
for any Borel subset $A$ of $M$.
The measure is called an {\em absolutely continuous conditionally invariant measure} (abbreviated 
a.c.c.i.m.)\ if its conditional distributions on unstable leaves are absolutely continuous
with respect to the Riemannian volume.

The quantity $\lambda = \mu (T^{-1}M)$ is called the eigenvalue of the measure and $-\log \lambda$
represents the exponential rate at which mass escapes from the system.
We call $\mu$ a {\em trivial} conditionally invariant measure if $\lambda = 0$.
From the point of view of physical observables, we are interested in conditionally invariant measures
whose escape rate indicates the rate at which the Riemannian volume escapes from the system.
For this reason, we will restrict our attention to the existence and properties of nontrivial absolutely
continuous conditionally invariant measures in this paper.

In this context, Pianigiani and Yorke \cite{pianigiani yorke} studied $C^2$ expanding maps which admit a 
finite Markov partition after the introduction of holes.  
Their work was extended to smooth Smale horseshoes by \v{C}encova (\cite{cencova 1}, \cite{cencova 2}) 
and used in \cite{collet 1} and \cite{collet 2} to obtain a natural invariant
measure on the fractal set of non-wandering points of the system.  
In \cite{lopes mark} and \cite{richardson}, the authors achieved similar results for open billiards 
satisfying a non-eclipsing condition.  
Recently, Chernov and Markarian (\cite{chernov mark1}, \cite{chernov mark2}) studied Anosov diffeomorphisms 
with holes which were elements of a finite Markov partition.  
In \cite{chernov mark t1} and \cite{chernov mark t2}, the Markov restriction on the holes 
was relaxed, but the results still used strongly the Markov partitions associated with Anosov diffeomorphisms.

In low-dimensional settings, efforts to drop the Markov requirements on both the map and the holes
have had some success for expanding maps of the interval.  
A spectral analysis of the transfer 
operator was performed in \cite{baladi keller} and the stability of the spectrum was established in 
\cite{keller liverani} for perturbations of expanding maps including small holes.  
More constructive techniques using bounded
variation and contraction mapping arguments have been used in \cite{chernov exp} and \cite{liverani maume} to 
prove the existence and properties of a.c.c.i.m.  
All these results assume that the potential associated with the transfer operator has bounded variation.

This brief survey highlights the classes of systems with holes which have been studied to date:
expanding maps in one dimension; and in higher dimensions, systems which admit finite Markov partitions.  
These systems are all uniformly hyperbolic. 
The purpose of this paper is to develop a method to study dynamical systems with holes which relies on 
neither finite Markov partitions nor uniform hyperbolicity.  
The Markov extension is a flexible tool which does precisely this.

The systematic application of Markov extensions represented schematically by tower maps is due to Young 
who has used this method to study a variety of closed dynamical systems including Axiom A diffeomorphisms, 
piecewise hyperbolic maps, 
H\'{e}non maps, logistic maps, and a class of scattering billiards, among others 
(\cite{young exp},\cite{young poly}, \cite{benedicks young}).  
Chernov has recently extended the method to study the statistical properties of other chaotic systems 
(\cite{chernov 1},\cite{chernov 2},\cite{chernov 3}).  
The flexibility of the Markov extension stems
from the fact that the dynamical system in question need not be uniformly hyperbolic.  
What matters is the average behavior of the map $T$ between returns to a reference set $\Lambda$.  
This is what allows the method to be applied to H\'{e}non maps and the logistic family.  
We briefly describe the main ideas of the construction.

Given a dynamical system $T: M \circlearrowleft$,  we choose a reference set $\Lambda \subset M$ and
run the system, waiting for a forward iterate of $\Lambda$ to make a ``good'' (Markov) return to $\Lambda$.
In the context of expanding systems, a good return is one which completely covers $\Lambda$.  
A stopping time is then declared on the piece of $\Lambda$ which makes the good return and we continue to iterate the
rest until another good return is made.  
In this way we generate a countable partition $\Z$ of $\Lambda$ and a stopping time $R$ which is constant on elements of $\Z$. 
Schematically, we consider $\Lambda$ as the
base $\Delta_0$ of a tower and each subsequent level $\Delta_l$ is identified with the part of $T^l(\Lambda)$
which has not made a good return to $\Lambda$ by time $l$.
The related system $(F, \Delta)$ is an extension of $T: \cup_{i\geq 0}T^i(\Lambda) \circlearrowleft$
for which there is a projection $\pi: \Delta \rightarrow \cup_{i\geq 0}T^i(\Lambda)$ which
satisfies $\pi \circ F = T \circ \pi$. 

There are three basic steps which are required for this method to yield results.
\begin{enumerate}
  \item[(1)] Given a dynamical system $T: M \circlearrowleft$, we construct a Markov extension
    $F : \Delta \circlearrowleft$;
  \item[(2)] we prove results about $(F, \Delta)$ using its simpler properties:  namely,
    controlled hyperbolicity and a countable Markov structure with a certain decay rate in the measure
    of the elements of the partition;
  \item[(3)] we pass these results back to the original system $(T, M)$.
\end{enumerate}
In Section~\ref{tower section}, we focus on (2) for a tower with holes.  
In Section~\ref{expanding section}, we apply the abstract model to expanding maps of the interval by proving (1) and (3).  
The author has also used Markov extensions to study logistic maps with holes.  
These results appear in a separate paper \cite{demers}.

\subsection{Summary of Main Results}
\label{summary of results}

Since important definitions and assumptions have not yet been introduced, here we state our results in general terms.  
Precise statements are made in Sections~\ref{tower section} and \ref{expanding section}.

Our principal result concerning Markov extensions is the following: 
\begin{enumerate}
        \item[(1)] Given a tower $(F, \Delta)$ with exponential decay in the levels and sufficiently 
        small holes, the tower map $F$ admits a conditionally invariant measure whose density is 
        H\"{o}lder continuous on elements of the Markov partition on $\Delta$.  If in addition $F$
        is transitive on components of $\Delta$, then the density is unique in the space of
        H\"{o}lder continuous functions and is bounded away from zero on $\Delta$.
\end{enumerate}
This result is stated precisely as Theorem~\ref{tower theorem} in Section~\ref{tower results}.
\\
\\
\noindent
{\bf Remark.} Because $\Delta$ is not compact, the exponential decay in the levels of the tower is 
essential to the discussion of meaningful conditionally invariant measures.  
Without it the escape rate of Lebesgue measure from $\Delta$ may not be exponential, thus making the 
conditionally invariant measure a poor indicator of the escape dynamics of the system.  
The exponential decay together with the transitivity assumption 
assures that the exponential escape rate of Lebesgue measure is well-defined.
\\
\\
In Section~\ref{expanding section}, we apply this result to expanding maps to obtain the following:
\begin{enumerate}
        \item[(2)] Let $T$ be a $C^{1+ \alpha}$ expanding map of the unit interval with sufficiently
        small holes.  Then $T$ admits an absolutely continuous conditionally invariant measure whose
        density is bounded above.  If $\alpha = 1$, then the density is of bounded variation.  If
        in addition $T$ satisfies a transitivity condition, the density is bounded away from zero on
        the complement of the hole, thus making the escape rate of Lebesgue measure equal to the
        escape rate of the a.c.c.i.m.  Moreover, the eigenvalue $\lambda$ satisfies the Lipschitz condition
        $1- \lambda \leq \mbox{const.}m(H)$.
        \item[(3)] If the measure of the holes shrinks to zero while the number of holes remains 
        bounded, then the a.c.c.i.m. converges weakly to the SRB measure of the expanding map without
        holes (here we assume the expanding map without holes has a unique a.c.i.m.).  If
        $\alpha = 1$, the convergence is in $L^1([0,1])$.
\end{enumerate}
These results are stated with full details as Theorems~\ref{expanding existence} and \ref{convergence}, 
respectively, in Section~\ref{expanding results}.
\\
\\
\noindent
{\bf Remark.}  The results of Theorem~\ref{expanding existence} are proved under more technical
assumptions in \cite{chernov exp} and \cite{liverani maume} for expanding maps having potentials of bounded variation.  
We do not use bounded variation arguments and only require that the potential be H\"{o}lder continuous.
Our only restriction on the hole is on its size, not its placement.

\section{An a.c.c.i.m.\ for a Tower with Holes}
\label{tower section}

We begin by establishing some basic notation and definitions regarding the tower map, following \cite{young exp}.
An important difference between the tower we shall describe here and that in \cite{young exp} is that we allow
multiple intervals in the base.  The returns are still Markov with respect to these intervals, but
need not cover the entire base.  Since we will apply these results in Section~\ref{expanding section} to
expanding maps of the interval, we formulate our definitions for one-dimensional
towers, but it is clear that our arguments can be extended to higher dimensional expanding towers as well.

\subsection{Setting and Assumptions}
\label{setting}

\subsubsection{Tower with multiple bases}
\label{multiple bases}

The base of the tower is an interval $\hat{\Delta}_0$ which is taken to be the union of finitely many
closed unit intervals whose interiors are pairwise disjoint, $\hat{\Delta}_0 = \bigcup_{i=1}^{N} \hdi_0$. 
We denote one-dimensional Lebesgue measure by $m$.
Let $\Z$ be a countable partition of $\hat{\Delta}_0$ each element of which is a subinterval of
one of the $\hdi_0$.   Given a return time function
$R:\hat{\Delta}_0 \rightarrow \mathbb{Z}^+$ which is constant on each
element of $\Z$, we define a tower $(\hat{\Delta}, \hat{F}, m)$ over 
$\hat{\Delta}_0$ by
\[
\hat{\Delta} := \{ (z,n) \in \hat{\Delta}_0 \times \mathbb{N} \; | \;
                    n < R(z) \}.
\]
We call the $l^{th}$ level of the tower 
$\hat{\Delta}_l := \hat{\Delta} | _{n=l}$ and $\hdi_l$ is the part of 
$\hat{\Delta}_l$ directly over $\hdi_0$.  Note that each level $\hat{\Delta}_l$ is simply
a union of disjoint intervals; it is not necessarily connected since some subintervals will
have returned to the base by time $l$. 
We let $\hdi = \bigcup_{l=0}^{\infty} \hdi_l$.

The action of $\hat{F}: \hat{\Delta} \rightarrow \hat{\Delta}$ is characterized by
\begin{eqnarray*}
        \F(z,l) = (z,l+1) & \! \! \! & \mbox{if $l+1 < R(z)$ and} \\
        \F^{R(z)}(\Z(z)) = \hdi_0 & \! \! \! & \mbox{for some $1 \leq i \leq N$}
\end{eqnarray*} 
where $\hat{F}^{R(z)} | _{\Z(z)}$ is continuous and one-to-one and $\Z(z)$ is the element
of $\Z$ containing $z$.  We adopt a slight abuse of notation by referring to $\F(x,l)$ as $\F(x)$
and $\Delta_l$ will be made clear by context.

As mentioned in the remark in Section~\ref{summary of results},
an essential requirement is an exponential decay in the measure of the levels of the tower.

\vspace{11 pt}
\noindent
(H1)  \hfill There exist $ A>0$ and $0< \theta < 1$ such that 
$\displaystyle m(\hat{\Delta}_l) \leq A \theta^l $ for $ l \geq 0$.
\hfill $\mbox{}$
\vspace{11 pt}

\subsubsection{Introduction of holes and regularity of $\mathbf{ \hat{F} }$ }
\label{regularity of F}

The hole $\h$ in $\hat{\Delta}$ is a union of open intervals $\hilj$.
We set $ \bigcup_j \hilj =: \h_l^{(i)} \subset \hdi_l$ and require that there be only finitely many $\hilj$
per level $l$ of the tower.  We define 
$\h_l = \bigcup_{i=1}^{N} \h_l^{(i)}$.  We also require that each 
$\hilj = \hat{F}^l(\omega)$, where $\omega$ is a union of elements of $\Z$, thus preserving the
Markov structure of the returns to $\hat{\Delta}_0$.  If 
$\hat{F}^l(\omega) = \hilj$, then the intervals on all levels of the tower 
directly above $\hilj$ are deleted since once $\hat{F}$ maps a point 
into $\h$, it disappears forever.  The interval $\omega$ does not return to $\hat{\Delta}_0$. 

Let $\Delta = \hat{\Delta} \backslash \h$ and 
$\Delta_l = \hat{\Delta}_l \backslash \h$ with analogous definitions for 
$\di$ and $\di_l$.  We assume the existence of a countable
Markov partition $\{ \dilj \}$ with $\bigcup_j \dilj = \di_l$ for each $i$ and $l$.  
Each $\dilj$ is an interval comprised of countably many elements of the form $\F^l(\omega)$, $\omega \in \Z$,
and $\F|_{\dilj}$ is one-to-one.

In applications,  $\{ \dilj \}$ will be dynamically defined during the construction
of the tower and its elements will be the maximal intervals which project onto the
iterated pieces of the reference set $\Lambda$ at time $l$.  Uniformly hyperbolic
systems, such as the expanding maps of Section~\ref{expanding section}, will typically
have finitely many $\dilj$ per level $l$.  We do not impose a finiteness requirement
here, however.

We denote by $\diljs$ those $\dilj$ whose image returns to the base, i.e. such that 
$\hat{F}(\dilj) = \bigcup_{k=k_1}^{k_2} \Delta^{(k)}_0$ for some $1 \leq k_1 \leq k_2 \leq N$
and set $\Delta^* = \bigcup \diljs$. 

The map $\hat{F}$ has the following properties with respect to the
partition:
\begin{enumerate}
   \item[(a)] $\hat{F}$ is $C^{1+\alpha}$ on each $\dilj$ for a fixed $\alpha >0$.
   \item[(b)] There exist $\gamma >1$ and $\beta > 0$ such that on $\diljs$, 
        $|\hat{F}'| \geq \gamma e^{\beta l}$.  Elsewhere, $|\hat{F}'|=1$.
   \item[(c)] {\em Bounded Distortion}.  There exists $C>0$
     such that for any $x,y \in \hdi_0$ 
     and $x', y' \in \dkljs$ such that 
     $\hat{F}(x')=x$ and $\hat{F}(y')=y$ we have
     \begin{equation}
         \left| \frac{\hat{F}'(x')}{\hat{F}'(y')} - 1 \right| 
         \leq C |x-y|^{\alpha}.             \label{eq:distortion}
     \end{equation}
\end{enumerate}
In applications, property (b) follows quite naturally once bounded distortion and the exponential decay of
the levels are proven.  In fact, $e^{\beta} > \frac{1}{\theta}$ and a variant of (b) is often used to prove (H1) in the
construction.  For uniformly hyperbolic systems it is immediate, but even for nonuniformly hyperbolic systems
such as logistic maps, it is easy to obtain once (c) and (H1) are established (see \cite{demers}).

We require an additional property of the map $\hat{F}$.  Taking
$C$, $\gamma$ and $\alpha$ as defined above in (a), (b) and (c),
we require:

\vspace{11 pt}
\noindent
(H2) \hfill  $\displaystyle \frac{1+C}{\gamma^{\alpha}} < 1.$
     \hfill $\mbox{}$
\vspace{11 pt}

This assumption says that the non-linearity of $\hat{F}$ should not
be large compared to its minimum expansion.  It is controlled by shrinking the size of the reference set
in the underlying dynamical system.
\\
\\
\noindent
{\bf Remark.}  If $\alpha = 1$, then
assumption (H2) is not necessary.  The justification for this is
contained in the remark at the end of Section~\ref{estimates}.
\\
\\
Let $F = \hat{F}|(\Delta \backslash \hat{F}^{-1}\h )$.  Note that $F^m$ is only defined on $\bigcap_{k=0}^m \F^{-k} \Delta$,
the set of points which has not fallen in the hole by time $m$.  We say $F$
is {\em transitive on components} if for all pairs $i,j$, $1 \leq i,j \leq N$,
there exists an $m$ such that $F^m( \di_0 \cap \bigcap_{k=0}^m \F^{-k} \Delta ) \supseteq \Delta^{(j)}_0$.  Note that
if $N=1$, then transitivity on components is automatic as long as the hole allows
at least one return to the base.

\subsubsection{Definition of a combined H\"{o}lder-$\mathbf{ L^{\infty} }$ Functional}

The last element of our setting is a
function space $X$ on $\Delta$ in which we will seek our conditionally
invariant density.  The Perron-Frobinius operator associated with $F$ acts on 
$L^1(\Delta)$ by
\[ \pf(x) = \sum_{y \in F^{-1}x} \frac{f(y)}{|F'(y)|}. \]
We define $\mathcal{P}_1f =  \pf / | \pf |_{L^1}$ and seek
a fixed point for the operator $\mathcal{P}_1$.  A fixed point for
$\mathcal{P}_1$ is a conditionally invariant density for $F$.

Choose $\xi > 0$ small enough that $e^{-\xi} > \max \{ \theta, e^{-\alpha \beta} \}$.
This is the only restriction we will put on $\xi$. 
Given $f \in L^1(\Delta)$, let $f^{(i)}_{l,j} = f|_{\dilj}$.  Let $|f|_{\infty}$ 
denote the $L^{\infty}$ norm of $f$ and define
\[ 
\|f^{(i)}_{l,j} \|_{\infty} = e^{-\xi l} |f^{(i)}_{l,j}|_{\infty}, 
\]
\[ 
\|f^{(i)}_{l,j} \|_r = \sup_{\stackrel{x,y \in \dilj}{f(x) \neq 0}} 
    \frac{|f(x)-f(y)|}{|x-y|^{\alpha} |f(x)| } e^{-\xi l}. 
\]
Then define
\[ 
\|f \| = \max \{ \|f\|_{\infty} , \|f \|_r  \} 
\]
where $\|f\|_{\infty} = \sup_{i,l,j} \|f^{(i)}_{l,j} \|_{\infty}$
and $\|f \|_r = \sup_{i,l,j} \|f^{(i)}_{l,j} \|_r$.
It should be noted that while $\| \cdot \|_{\infty}$ is a norm, $\| \cdot \|_r$
is not; however, $\| \cdot \|_r$ satisfies a convex-like inequality on a subset of $X$ defined
in Section~\ref{convex subset}.  We set 
\[
X = \{ f:\Delta \rightarrow \mathbb{R} \; | \; \| f \| < \infty \}.
\]

\subsubsection{Condition on the Holes}
\label{condition on holes}

We formulate a single condition involving the measure of the holes which
guarantees the existence of an a.c.c.i.m.\ in X.

Let $a:= \max \{ e^{-\xi},  
\frac{1+C}{\gamma^{\alpha}} \}$ and
$b:= 1+C $.
Note that $a<1$ by assumption (H2).  The required condition on the 
holes is:

\vspace{11 pt}
\noindent
(H3) \hfill  $\displaystyle Dm\h_0 + \sum_{l \geq 1} e^{\xi (l-1)} m\h_l
             \leq \frac{(1-a)^2}{4b}$
     \hfill $\mbox{}$
\vspace{11 pt}

\noindent
where $D = (1+C) \sum_{\dkljs} e^{\xi l}m \dkljs$.  In typical
applications of the tower constructed from a dynamical system, 
there will be no holes in the base of the tower.  In this case
$m\h_0 = 0$, the proof simplifies somewhat and assumption (H3) reduces to

\vspace{11 pt}
\noindent
(H$3'$) \hfill  $\displaystyle \sum_{l \geq 1} e^{\xi (l-1)} m\h_l
             \leq \frac{(1-a)^2}{b}.$
     \hfill $\mbox{}$

\subsubsection{Statement of results}
\label{tower results}

We are now ready to state our main result of this section.
\begin{theorem}
\label{tower theorem}
Given a tower with holes $(\Delta, F, m)$ with properties (a)--(c)
and under assumptions (H1)--(H3),
there exists a probability density $\varphi$ in $X$ such
that $\mathcal{P}_1 \varphi = \varphi$.  If in addition $F$ is
transitive on components, then $\varphi$ is the unique nontrivial conditionally invariant density
in $X$ and $\varphi$ is bounded away from zero on $\Delta$.
\end{theorem}

\subsection{Proof of Theorem \ref{tower theorem}}

Our proof takes the following steps:
\begin{enumerate}
  \item define a convex, compact subset $X_M$ of $X$; 
  \item derive Lasota-Yorke type inequalities for the operator $\mathcal{P}$;
  \item formulate the condition on $\h$ which guarantees the invariance of
    $X_M$ under $\mathcal{P}_1$; 
  \item use transitivity to prove the uniqueness of the invariant density in $X$.
\end{enumerate}

\subsubsection{Properties of the space $\mathbf{X}$ and a convex, compact subset}
\label{convex subset}

First note that $X \subset L^1(\Delta)$.  For
if we let $|\cdot|_1$ denote the norm in $L^1(\Delta)$, then
for $f \in X$ we have
\[ 
|f|_1 = \sum_{i,l,j} \int_{\dilj} |f| dm 
    \leq \sum_{i,l,j} |f^{(i)}_{l,j}|_{\infty} m(\dilj)
    \leq \|f\|_{\infty} \sum_{i,l,j} m(\dilj) e^{\xi l}, 
\]
which is finite by the assumption on $\xi$.  We also record for future use that if
$f \in X$, then for $x,y \in \dilj$
\[ 
\|f^{(i)}_{l,j} \|_r \geq \frac{|f(x)-f(y)|}{|x-y|^{\alpha} |f(x)|} e^{-\xi l}
      = \frac{ \left| \frac{f(y)}{f(x)} - 1 \right|}{|x-y|^{\alpha}} e^{-\xi l} .  
\]
This in turn yields
\begin{equation}
  \sup_{\dilj} |f| \leq (1 + \| f \|_r e^{\xi l} ) \inf_{\dilj} |f|  \label{eq:supinf}
\end{equation}
so that $f$ is either identically zero on $\dilj$ or $f$ is bounded 
away from zero on $\dilj$.

Now let 
\[ 
X_M = \{ f \in X \; | \; \| f \| \leq M, f \geq 0, \int_{\Delta} f dm =1  \} 
\]
where $M>0$ has yet to be determined.

\begin{proposition}
$X_M$ is a convex, compact subset of $L^1(\Delta)$.
\end{proposition}
{\em Proof}.  We first show that $X_M$ is convex.  Let $f,g \in X_M$.  Then
$\| sf + (1-s)g \|_{\infty} \leq s \| f \|_{\infty} + (1-s) \| g \|_{\infty}
\leq M$ for $s \in [0,1]$ since $\| \cdot \|_{\infty}$ is a norm.

To show that $\| \cdot \|_r$ satisfies a similar convexity property, we use the
fact that for any positive numbers $a_1, a_2, b_1, b_2$,
\begin{equation}
\label{eq:convex hull}
\frac{a_1 + a_2}{b_1 + b_2} \leq \max \left\{ \frac{a_1}{b_1}, \frac{a_2}{b_2} \right\}.
\end{equation}
If $f \equiv 0$ on $\dilj$, then clearly $\|(sf + (1-s)g)^{(i)}_{l,j} \|_r = \| g^{(i)}_{l,j} \|_r \leq M$.
So let us assume that both $f$ and $g$ are positive on $\dilj$.  Then
\begin{eqnarray*}
\|(sf + (1-s)g)^{(i)}_{l,j} \|_r 
  & \leq & \sup_{x,y \in \dilj}\frac{ s|f(x)-f(y)| + (1-s)|g(x) - g(y)| } { |x-y|^{\alpha} (sf(x) + (1-s)g(x)) } e^{-\xi l} \\
  & \leq & \max \{ \| f^{(i)}_{l,j} \|_r , \| g^{(i)}_{l,j} \|_r \} \; \leq \; M
\end{eqnarray*}
where we have used equation~(\ref{eq:convex hull}) for the second inequality.

Thus $sf + (1-s)g \in X_M$ whenever $f,g \in X_M$, so $X_M$ is convex.  We now show that
$X_M$ is compact as well.

Let $\{ f_n \}_{n=0}^{\infty}$ be a sequence in $X_M$.
For fixed $i$, $l$ and $j$, 
\[  
\sup_{x,y \in \dilj} \frac{|f_n(x) - f_n(y)|}{|x-y|^{\alpha}}
       \leq \| f_n \|_r \|f_n \|_{\infty} e^{2 \xi l}.  
\]
This makes $\{ f_n|_{\dilj} \}$ a uniformly Holder continuous and thus
equicontinuous family of functions.  It is also bounded in $L^{\infty}$.
By the Ascoli-Arzela theorem, there exists a subsequence
which converges pointwise uniformly on $\dilj$ to a function $f_{*l,j}^{(i)}$.
By the uniformity of the convergence, $f_{*l,j}^{(i)}$ has the same 
$\| \cdot \|_r$ bound as the $f_n|_{\dilj}$.  
It also has the same bound in $L^{\infty}$.

Using Cantor diagonalization, we obtain a subsequence which converges
pointwise on all of $\Delta$. Call this sequence $\{ f_{n_k} \}$ and
the limit $f_*$.  Since $|f_{n_k}|$ is dominated by $Me^{\xi l}$ on each $\Delta_l$,
we conclude that $f_{n_k} \rightarrow f_*$ in $L^1(\Delta)$ as well.  
\\
\\
Note that our choice of $\xi$ already limits the size of $\h$ by
limiting the eigenvalue of the invariant measure.  Suppose we find
a nontrivial conditionally invariant density $\varphi \in X_M$ and let $\varphi dm = d\mu$.  
If $\lambda := \mu (F^{-1}\Delta) <1$, then we must have
$\varphi |_{\Delta_{l+1}} = \frac{1}{\lambda} \varphi |_{F^{-1}(\Delta_{l+1})}$ by the
conditional invariance property.  Since $\varphi \in X_M$ and in particular
$\| \varphi \|_{\infty}$ is finite, we must have 
$\lambda \geq e^{-\xi}$.

\subsubsection{Estimates on $\mathbf{\|\pf\|}$}
\label{estimates}

Our goal is to show that there exists a choice of
$M$ such that $\mathcal{P}_1$ takes $X_M$ into itself.  We do this 
by first deriving Lasota-Yorke type inequalities for the operator $\mathcal{P}$.
Note that by definition, $\pf(x)  = f(F^{-1}x)$ for $x \in \Delta_l$, $l \geq 1$.

For $\dilj$, $l \geq 1$, we have the estimate,
\[
\begin{array}{ccccl}
  \displaystyle \| \pf^{(i)}_{l,j} \|_{\infty} & := 
        & \displaystyle \sup_{x \in \dilj} | \pf (x)|e^{-\xi l}
        & = & \displaystyle \sup_{x \in \dilj} |f(F^{-1}x)| e^{-\xi l} \\
    & = &  \displaystyle 
        \sup_{y \in F^{-1}\dilj} |f(y)|e^{-\xi (l-1)}e^{-\xi}
        & \leq & \displaystyle \| f \|_{\infty}e^{-\xi}.
\end{array}
\]
And similarly,
\[
\begin{array}{ccccl}
  \| \pf^{(i)}_{l,j} \|_r & = & \displaystyle \sup_{x,y \in \dilj} 
            \frac{|f(F^{-1}x)-f(F^{-1}y)|}{|x-y|^{\alpha} |f(F^{-1}x)|} e^{-\xi l} & & \\
    & = & \displaystyle \sup_{w,z \in F^{-1}\dilj}
            \frac{|f(w)-f(z)|}{|w-z|^{\alpha} |f(w)|}  e^{-\xi l}  
    & \leq & \displaystyle \| f \|_r e^{-\xi} .
\end{array}
\]

For $\di_0$ the estimates are more involved.  First note that 
for $x,y \in \dkljs$ with $F(x),F(y) \in \di_0$, equation~(\ref{eq:distortion}) yields
\begin{equation}
  \left| \frac{1}{|F'(x)|} - \frac{1}{|F'(y)|} \right| \leq
         C |F(x)-F(y)|^{\alpha}\frac{1}{|F'(y)|}          \label{eq:difference}
\end{equation}
so that
\begin{equation}
  \frac{1}{|F'(x)|} \leq (1+C) m\dkljs             \label{eq:mdlj}
\end{equation}
and also
\begin{equation}
  \left| \frac{F'(y)}{F'(x)} \right| 
        \leq 1 + C .                                   \label{eq:ratio}
\end{equation}

\noindent
{\bf Estimating $\mathbf{\| \pf^{(i)}_{0,j} \|_{\infty} }$.} \\
For $x \in \di_0$ we have $\pf (x) = \invx \frac{f(y)}{|F'(y)|}$.  For
a fixed $\dkljs$, let $\bar{a} \in \dkljs$ be such that 
$f(\bar{a}) = \frac{1}{m\dkljs}\int_{\dkljs}fdm$.  Then
\begin{equation}
  |f^{(k)}_{l,j}|_{\infty} \leq \frac{1}{m\dkljs} \int_{\dkljs} f dm
          + \sup_{x \in \dkljs} |f(x) - f(\bar{a})| .   \label{eq:linftysplit}
\end{equation}
We use (\ref{eq:mdlj}) to estimate the first term,
\begin{equation}
   \frac{1}{|F'(y)| m\dkljs} \int_{\dkljs} f dm \leq
          (1+C) \int_{\dkljs} f dm.       \label{eq:inftyfirst}
\end{equation}
We estimate the second term by,
\begin{eqnarray*}
  \sup_{x \in \dkljs} |f(x)-f(\bar{a})| & = & 
  \sup_{x \in \dkljs} \frac{|f(x)-f(\bar{a})|}{|x-\bar{a}|^{\alpha} 
                            f(\bar{a}) } 
          |x-\bar{a}|^{\alpha} f(\bar{a}) \\
  & \leq & \| f^{(k)}_{l,j} \|_r e^{\xi l} (m \dkljs )^{\alpha}
          \frac{1}{m\dkljs}\int_{\dkljs}fdm.
\end{eqnarray*}
We then use (\ref{eq:mdlj}) and property (b) to split up $|F'|$,
\[
\frac{1}{|F'(y)|} = 
   \frac{1}{|F'(y)|^{1-\alpha}} \cdot \frac{1}{|F'(y)|^{\alpha}}
   \leq (m\dkljs)^{1-\alpha} (1+C)^{1-\alpha}
   \frac{1}{\gamma^{\alpha}} e^{-\alpha \beta l}
\]
so that
\begin{equation}
  \frac{\sup_{x \in \dkljs} |f(x)-f(\bar{a})|}{|F'(y)|}
  \leq \| f^{(k)}_{l,j} \|_r  \frac{(1+C)^{1-\alpha} }{\gamma^{\alpha}} 
       \int_{\dkljs}fdm.      \label{eq:inftysecond}
\end{equation}
where we have used the fact that $e^{\xi} \leq e^{\alpha \beta}$.
Combining (\ref{eq:linftysplit}), (\ref{eq:inftyfirst}) and 
(\ref{eq:inftysecond}) we obtain:
\begin{eqnarray}
  \|\pf^{(i)}_{0,j} \|_{\infty} & = & |\pf^{(i)}_{0,j}|_{\infty} = 
                \sup_{x \in \di_{0,j}} \invx \frac{f(y)}{|F'(y)|}  \nonumber \\
       & \leq & \sum_{\dkljs} (1+C)
                \int_{\dkljs}fdm +  \| f^{(k)}_{l,j} \|_r 
                \frac{ (1+C)^{1-\alpha}}{\gamma^{\alpha}}
                \int_{\dkljs}fdm                                  \nonumber \\
       & \leq & (1+C) \int_{\Delta^*} f dm + \|f\|_r
                \frac{ (1+C)^{1-\alpha}}{\gamma^{\alpha}}
                \int_{\Delta^*} fdm                        \label{eq:integral bound} \\
       & \leq & \frac{ (1+C)^{1-\alpha}}{\gamma^{\alpha}}\|f\|_r + 1+C .     \nonumber
\end{eqnarray}

\noindent
{\bf Estimating $\mathbf{ \|\pf^{(i)}_{0,j} \|_r }$.} \\
Let $x,z \in \di_{0,j}$.
\begin{eqnarray*}
  | \pf(x) - \pf(z) | & \leq & \invxz \left|  \frac{f(y)}{|F'(y)|} - 
             \frac{f(w)}{|F'(w)|} \right|  \\
         & \leq & \invxz \left| \frac{f(y)}{|F'(y)|} - 
             \frac{f(y)}{|F'(w)|} \right| + \left| \frac{f(y)}{|F'(w)|} -
             \frac{f(w)}{|F'(w)|} \right|  \\
         & =: & \sum_{\dkljs} A^k_{l,j} + B^k_{l,j}.
\end{eqnarray*}

We use the following extension of equation~(\ref{eq:convex hull}):  given two convergent series, $\sum a_i$
and $\sum b_i$, $a_i \geq 0$, $b_i > 0$, then
\begin{equation}
  \frac{\sum a_i}{\sum b_i} \leq \sup_i \frac{a_i}{b_i}. \label{eq:fact}
\end{equation}
Note that if $f(y)=0$, then $f(w)=0$ by (\ref{eq:supinf}) since
$y$ and $w$ are in the same $\dklj$.
So dropping terms where $f(y)=0$, we may use (\ref{eq:fact}) to estimate
\begin{eqnarray*}
  \frac{|\pf(x) - \pf(z)|}{|x-z|^{\alpha} |\pf(x)|} & \leq &
         \frac{ \sum_{\dkljs} A^k_{l,j} + \sum_{\dkljs} B^k_{l,j} }
              { |x-z|^{\alpha} \invx \frac{f(y)}{|F'(y)|} }  \\
       & \leq & \sup_{\stackrel{\dkljs}{F(y)=x}} 
              \frac{A^k_{l,j}}{|x-z|^{\alpha} \frac{f(y)}{|F'(y)|}}
              + \sup_{\stackrel{\dkljs}{F(y)=x}} 
              \frac{B^k_{l,j}}{|x-z|^{\alpha} \frac{f(y)}{|F'(y)|}}
\end{eqnarray*}
For the first term,
\[
\frac{A^k_{l,j}}{|x-z|^{\alpha} \frac{f(y)}{|F'(y)|}} 
   = \frac{1}{|x-z|^{\alpha}} \left| \frac{F'(y)}{F'(w)} - 1 \right|
   \leq C
\]
where we have used (\ref{eq:distortion}) for the inequality.
For the second term,
\[
\frac{B^k_{l,j}}{|x-z|^{\alpha} \frac{f(y)}{|F'(y)|}}
   = \frac{|f(y) - f(w)|}{|x-z|^{\alpha} f(y)} \cdot
     \frac{|F'(y)|}{|F'(w)|} 
   \leq \| f^{(k)}_{l,j} \|_r \frac{1}{\gamma^{\alpha}} (1 + C)
\]
where we have used the fact that
$|x-z| \geq \left( \inf_{\dkljs} |F'| \right) |y-w|$ in addition to equation (\ref{eq:ratio}) and property (b),
Combining the estimates for $A^k_{l,j}$ and $B^k_{l,j}$ we see that
\[
\| \pf^{(k)}_{0,j} \|_r \leq C +
        \| f \|_r \frac{1}{\gamma^{\alpha}} (1 + C) .
\]

Recalling that
$a = \max \{ e^{-\xi}, \frac{1+C}{\gamma^{\alpha}} \}$ and
$b = 1 + C$
we have shown that
\begin{equation}
\label{r estimates}
\| \pf \|_r \leq a \| f \|_r + b .
\end{equation}
We also have
\begin{equation}
\label{infinity estimates}
\| \pf^{(i)}_{0,j} \|_{\infty} \leq  a \| f \|_r + b \; \; \; \mbox{and} \; \; \;
          \| \pf^{(i)}_{l,j} \|_{\infty} \leq a \|f \|_{\infty} \; \; \mbox{for  $l \geq 1$}
\end{equation}
where $a<1$ by assumption (H2).
\\
\\
\noindent
{\bf Remark.}  If $\alpha = 1$, i.e. $F$ is $C^2$ on each $\dilj$, 
then we do not need assumption (H2).  We define 
$ \| f^{(i)}_{l,j}\|_r = e^{-\xi l} \sup_{x \in \dilj} \left\{ \frac{|f'(x)|}{|f(x)|} \;
                   | \; f(x) \neq 0 \right\}$.
Using this norm, we note that on $\dilj$,
\[ 
|f(x) - f(\bar{a})| = | \int_{\bar{a}}^x f' dm | 
          \leq \|f^{(i)}_{l,j} \|_r e^{\xi l} \int_{\dilj} f dm.
\]
Using this estimate, equation (\ref{eq:inftysecond}) becomes
\[
  \frac{\sup_{x \in \dkljs} |f(x)-f(\bar{a})|}{|F'(y)|}
  \leq \| f^{(k)}_{l,j} \|_r \frac{1}{\gamma} \int_{\dkljs}fdm.  
\]
So $\|\pf^{(i)}_{0,j}\|_{\infty} \leq \frac{1}{\gamma} \| f \|_r + C + 1$.

The estimate on $\| \pf^{(i)}_{0,j} \|_r$ changes slightly as well.  If we
denote by $\psi^{(k)}_{l,j}$ the inverse of $F|_{\dklj}$, then for 
$x \in \di_0$ we have
\begin{eqnarray*}
\frac{ |\pf'(x)|}{\pf(x)} & = &\frac{ \sum_{\dkljs} 
    |(f \circ \psi^{(k)}_{l,j}(x) \cdot |(\psi^{(k)}_{l,j})'(x)|)'| }
    { \sum_{\dkljs} f \circ \psi^{(k)}_{l,j}(x) \cdot |(\psi^{(k)}_{l,j})'(x)| } \\
  & \leq & \sup_{\dkljs} 
      \frac{ |f'| \circ \psi^{(k)}_{l,j}(x) \cdot |(\psi^{(k)}_{l,j})'(x)|^2 } 
           { f \circ \psi^{(k)}_{l,j}(x) \cdot |(\psi^{(k)}_{l,j})'(x)| }
      + \sup_{\dkljs} 
      \frac{ f \circ \psi^{(k)}_{l,j}(x) \cdot |(\psi^{(k)}_{l,j})''(x)| } 
           { f \circ \psi^{(k)}_{l,j}(x) \cdot |(\psi^{(k)}_{l,j})'(x)| }  \\
  & \leq & \frac{1}{\gamma} \| f \|_r + C
\end{eqnarray*}
where we have used (\ref{eq:fact}) in the second step
and (\ref{eq:difference}) in the second term of the third.
In this case, $a := \max \{ e^{-\xi}, \frac{1}{\gamma} \}$ which
is automatically less than one.

\subsubsection{Choice of $\mathbf{M}$ and condition on $\mathbf{\h}$}
\label{choice of M}

Choose $M = \frac{2b}{1-a}$.  We derive a condition on the size of the 
hole depending on $a, b$  and $\xi$ which guarantees that 
$\mathcal{P}_1$ will map $X_M$ into itself.  This will be precisely
assumption (H3) introduced in Section~\ref{condition on holes}.

Since $\| \cdot \|_r$ is invariant under constant multiples of a 
given function $f$, we have $\| \mathcal{P}_1f \|_r = \| \pf \|_r$.
Since $b= \frac{1-a}{2} M$, using equation~(\ref{r estimates}) we have
\[
\| \textstyle \mathcal{P}_1f \|_r = \| \pf \|_r \leq aM + b
     = aM + \frac{1-a}{2} M < M.
\]

The $\| \cdot \|_{\infty}$ norm is not invariant under normalization,
however, so we need to introduce a restriction on the size of the holes. 
First note that for $f \in X_M$
\begin{eqnarray*}
  |\pf|_1 & = & \int_{\Delta \backslash (\hat{F}^{-1}\h)} fdm 
                 = 1 - \int_{\hat{F}^{-1}(\h) \cap \Delta} fdm \\
      & = & 1 - \sum_{l \geq 1} \int_{\hat{F}^{-1}\h_l} fdm - 
           \int_{\hat{F}^{-1}(\h_0) \cap \Delta} fdm \\
      & \geq & 1 - \sum_{l \geq 1} |f_{l-1}|_{\infty} m(\hat{F}^{-1}\h_l)
           - \sum_{\dkljs} \int_{\hat{F}^{-1}(\h_0) \cap \dkljs} fdm \\
      & \geq & 1-\sum_{l \geq 1} \| f_{l-1} \|_{\infty} e^{\xi(l-1)}m\h_l
           - \sum_{\dkljs} \| f_{l,j} \|_{\infty} e^{\xi l} m\h_0
             m\dkljs (1+C)  \\
      & \geq & 1 - \|f \| \left(  \sum_{l \geq 1}e^{\xi(l-1)}m\h_l
           + Dm\h_0 \right).
\end{eqnarray*}
where $D = (1+C) \sum_{\dkljs} e^{\xi l} m\dkljs$ as before.
To simplify this expression, let $q:= Dm\h_0 + \sum_{l \geq 1}e^{\xi(l-1)}m\h_l$ and note that
$q \rightarrow 0$ as $m\h \rightarrow 0$.  We require that $\h$ be small enough
so that 
\[
q \leq \frac{1-a}{2M} = \frac{(1-a)^2}{4b}
\]
which is precisely assumption (H3).
Using this bound on $\h$, we conclude that
\[
\begin{array}{rcccl}
  \displaystyle \| \mathcal{P}_1f \|_{\infty} & \leq & \displaystyle \frac{ a \| f \| + b }{|\pf|_1}
             & \leq & \displaystyle \frac{ a \| f \| + b}{1 - q \|f\|} \vspace{4 pt} \\
             & \leq & \displaystyle \frac{aM + \frac{1-a}{2}M}{1 - \frac{1-a}{2M} M} 
             & = & M.
\end{array}
\]
Note that this estimate implies that $\lambda \geq \frac{1+a}{2}$.

In the case when there are no holes in the base, we may use the simpler assumption
(H$3'$).  Now choose $M = \frac{b}{1-a}$.  Repeating our previous argument, we have
$\| \mathcal{P}_1f \|_r \leq M$ when $\|f\|_r \leq M$.

For the $\| \cdot \|_{\infty}$ norm, note that equation (\ref{eq:integral bound})
yields 
\[
\| \pf_{0,j}^{(i)} \|_{\infty} \leq (a \|f \|_r + b)\int_{\Delta \backslash (\hat{F}^{-1}\h)} fdm
\]
since $\Delta^* \subseteq \Delta \backslash (\hat{F}^{-1}\h)$ when there are no holes
in the base.  But  $|\pf|_1 = \int_{\Delta \backslash (\hat{F}^{-1}\h)} fdm$ so we have
$\| \mathcal{P}_1f_{0,j}^{(i)} \|_{\infty} \leq aM+b \leq M$.

All that remains is the estimate on $\| \mathcal{P}_1f^{(i)}_{l,j} \|_{\infty}$
for $l \geq 1$.
\[
\begin{array}{rcccl}
\displaystyle \| \mathcal{P}_1f^{(i)}_{l,j} \|_{\infty} & \leq & \displaystyle \frac{aM}{|\pf|_1}
                    & \leq & \displaystyle \frac{aM}{1-qM}.
\end{array}
\]
This is less than or equal to $M$ if $q \leq \frac{1-a}{M} = \frac{(1-a)^2}{b}$,
which is precisely assumption (H$3'$).

Since $\mathcal{P}_1$ takes $X_M$ into itself and $X_M$ is a convex,
compact subset of $L^1(\Delta)$, we can apply the Schauder-Tychonov
Theorem to conclude that $\mathcal{P}_1$ has a fixed point 
$\varphi \in X_M$ once we demonstrate that $\mathcal{P}_1$ is a continuous map on $X_M$.  
This follows from the estimate
\begin{eqnarray*}
|\mathcal{P}_1 f - \mathcal{P}_1 g |_1 & \leq & \left| \frac{\mathcal{P}f}{|\mathcal{P}f|_1} 
        - \frac{\mathcal{P}g}{|\mathcal{P}f|_1} \right|_1 +
        \left| \frac{\mathcal{P}g}{|\mathcal{P}f|_1} 
        - \frac{\mathcal{P}g}{|\mathcal{P}g|_1} \right|_1  \\
        & \leq & \frac{|f-g|_1}{|\mathcal{P}f|_1} + \frac{\left| |\mathcal{P}g|_1 - |\mathcal{P}f|_1 \right|}
                                                                        {|\mathcal{P}f|_1}  \\
        & \leq & \frac{|f-g|_1}{|\mathcal{P}f|_1} + \frac{|\mathcal{P}f-\mathcal{P}g|_1}{|\mathcal{P}f|_1} 
                  \leq 2 \frac{|f-g|_1}{|\mathcal{P}f|_1} 
\end{eqnarray*}
together with the fact that $|\mathcal{P}f|_1 \geq 1 - qM$ for $f \in X_M$.
Setting $d\mu = \varphi dm$, we have a measure 
that is conditionally invariant with respect to $F$ and absolutely 
continuous with respect to Lebesgue measure.  The eigenvalue
$\lambda = |\mathcal{P} \varphi|_1$ satisfies $\lambda \geq 1- qM$
by our estimates above and so $\lambda \rightarrow 1$ as $m\h \rightarrow 0$.

\subsubsection{Uniqueness of the conditionally invariant measure}

We now prove that $\varphi$ is the unique nontrivial conditionally invariant density 
in the function space $X$ if $F$ is transitive on components.
Recall from (\ref{eq:supinf}) that if $f \in X$, then on each $\dilj$,
$\sup_{\dilj} |f| \leq (1 + \|f\|e^{\xi l} ) \inf_{\dilj} |f|$, so that $f^{(i)}_{l,j}$
is either identically zero or bounded away from zero.  For a 
conditionally invariant density $\varphi$ with eigenvalue
$0 < \lambda \leq 1$, if $\varphi \equiv 0$ on 
some $\di_{0,j}$,
then $\varphi \equiv 0$ on $\Delta$ by transitivity. And since if $\F(\dkljs)$ intersects $\di_0$, it
crosses it completely, we must have 
$\inf_{\Delta_0}\varphi = \delta$ for some $\delta >0$.  
But by the conditional invariance
property, $\varphi |_{\Delta_{l+1}} = \frac{1}{\lambda} \varphi|_{F^{-1}(\Delta_{l+1})}$ for $l \geq 0$.
This makes $\inf_{\Delta} \varphi = \inf_{\Delta_0} \varphi$.

\begin{proposition}
\label{regularity}
If $\varphi \in X$ is a nontrivial conditionally invariant density for $F$, then $\varphi \in X_M$.
\end{proposition}

\noindent
{\em Proof.}  
From equation~(\ref{r estimates}), we have $\| \mathcal{P}_1 \varphi \|_r \leq a \| \varphi \|_r + b$
which implies
\[
\|\varphi \|_r = \| \mathcal{P}_1^n \varphi  \|_r \leq a^n \| \varphi \|_r + \frac{b}{1-a}
\; \; \mbox{for all $n$}.
\]
Thus $\| \varphi \|_r \leq \frac{b}{1-a}$.

For the $\| \cdot \|_{\infty}$ norm, the conditional invariance property implies that
\[
\| \varphi^{(i)}_{l,j} \|_{\infty} = e^{-\xi l} |\varphi^{(i)}_{l,j}|_{\infty} 
 \leq e^{-\xi l} \frac{1}{\lambda^l} | \varphi|_{\Delta_0} |_{\infty} \leq \| \varphi |_{\Delta_0} \|_{\infty}
\]
since $\lambda \geq e^{-\xi}$.

From equation~(\ref{infinity estimates}), we have 
$\| \mathcal{P} \varphi_{0,j}^{(i)} \|_{\infty} \leq  a \| \varphi \|_r + b$.

If there are holes in the base, then $\lambda \geq \frac{1+a}{2}$ so
\[
 \| \mathcal{P}_1 \varphi_{0,j}^{(i)} \|_{\infty} =  \frac{\| \mathcal{P} \varphi_{0,j}^{(i)} \|_{\infty}}{\lambda}
      \leq \frac{ a \frac{b}{1-a} + b}{ \lambda} \leq \frac{2b}{1-a} .
\]

If there are no holes in the base, then equation~(\ref{eq:integral bound}) yields
\[
 \| \mathcal{P}_1 \varphi_{0,j}^{(i)} \|_{\infty} 
    \leq \frac{( a \| \varphi \|_r + b)}{\lambda} \int_{\Delta \backslash \hat{F}^{-1}\h } \varphi dm
    \leq a \frac{b}{1-a} + b = \frac{b}{1-a} .
\]

So whether there are holes in the base and we define $M = \frac{2b}{1-a}$ or there are
no holes in the base and we define $M = \frac{b}{1-a}$, we have $\varphi \in X_M$.
\hfill $\Box$   \\

Now assume that there are two different nontrivial conditionally invariant densities
$\vone$ and $\vtwo$ in $X$.  Suppose they both have the same eigenvalue
$\lambda$. For $s \in \mathbb{R}$, let $\varphi_s = s \vone + (1-s) \vtwo$.  
Then $\varphi_s$ is a conditionally invariant density for each $s$ as 
long as $\varphi_s > 0$.  This is an open condition since $\inf_{\Delta} \varphi_s = \inf_{\Delta_0} \varphi_s$
and will certainly be true for $s \in [0,1]$ since $X_M$ is convex.
Let $s_0$ be the first $s > 1$ such that $\inf_{\Delta} \varphi_{s_0} = 0$.
For each $s<s_0$, $\varphi_s \in X_M$ by Proposition~\ref{regularity} so that $\varphi_{s_0} \in X_M$
as well since $X_M$ is closed.  But then we must have  $\varphi_{s_0}$
identically equal to zero which is impossible since $|\varphi_{s_0}|_1 = 1$. 

Now suppose that $\vone$ and $\vtwo$ have different eigenvalues
$\lone > \ltwo$.  Since $\vtwo$ is bounded away from zero and $\vone$
is bounded above on $\Delta_0$, there exists a constant $L$ such that
$L \vtwo > \vone$ on $\Delta_0$. We also know that 
$\vone |_{\Delta_{l+1}} = \frac{1}{\lone} \vone|_{F^{-1}(\Delta_{l+1})}$ by the 
conditional invariance property and similarly for $\vtwo$.  Since
$\frac{1}{\ltwo} > \frac{1}{\lone}$, then $L \vtwo > \vone$ on all of
$\Delta$.  So by the positivity of the operator $\mathcal{P}$,
we must have $\mathcal{P}^n(L \vtwo) > \mathcal{P}^n\vone$ for each
$n \geq 0$.  But this implies that $\ltwo^n L \vtwo > \lone^n \vone$
for each $n$ which is impossible since $\lone > \ltwo$.

\section{Expanding Maps with Holes: An Application}
\label{expanding section}

We describe an application of the tower model with holes to a class of open chaotic dynamical systems.
By constructing a tower with multiple bases and applying the results of the previous section, we study 
the existence and properties of an absolutely continuous conditionally invariant measure 
for $C^{1+\alpha}$ expanding maps of the interval with holes.

\subsection{Setting and Statement of Results}
\label{expanding setting}

\subsubsection{Properties of the Expanding Map}
\label{properties}

Let $\T$ be an expanding map of the unit interval, $\I=[0,1]$.  Denote by $\hat{I}_j$, 
$1 \leq j \leq \hat{K}$, the intervals of monotonicity for $\T$.  We assume that
\begin{enumerate}
   \item[(a)] $\T$ is $C^{1+ \alpha}$ on each $\hat{I}_j$.
   \item[(b)] $|\T'| \geq \mu >2$.
\end{enumerate}
Property (a) implies that there exists $\hat{C}>0$ such that 
$| \T'(x) - \T'(y)| \leq \hat{C} |x-y|\alf $ whenever $x,y$ are in
the same $\hat{I}_j$.  We let $m$ denote Lebesgue measure on the tower and on $\I$ interchangeably.

Now let $\hat{I}_j^n$ be the intervals of monotonicity for $\T^n$.  Let
$\psi_{n,j}$ denote the inverse of $\T^n$ acting on $\hat{I}_j^n$.  Then 
$|\psi'_{n,j} | \leq 1/ \mu^n $ so that
\[
|x-y| \leq \frac{1}{\mu^n} |\T^n(x) - \T^n(y)|
\]
whenever $x,y$ are in the same $\hat{I}^n_j$.  This gives rise to the following familiar fact.

\begin{lemma}
\label{lem:distortion}
{\em (Distortion Bounds).}
Let $x,y \in \hat{I}^n_j$.  Then
\[
\left| \frac{(\T^n)'(x')}{(\T^n)'(y')} - 1 \right| \leq \C |\T^n(x)-\T^n(y)|\alf
\]
where $\C = \exp(\frac{\hat{C}}{\mu(\mu\alf - 1)}) -1$.
\end{lemma}

\noindent
{\em Proof}.  
\begin{eqnarray*}
\displaystyle \log \left| \frac{(\T^n)'(x)}{(\T^n)'(y)} \right|
      & \leq &  \sum_{i=0}^{n-1} | \log |\T'(\T^ix)| - \log|\T'(\T^iy)| |  \\
      & \leq &  \sum_{i=0}^{n-1} \frac{1}{\mu} |\T'(\T^ix) - \T'(\T^iy)|   \\
      & \leq &  \sum_{i=0}^{n-1} \frac{\hat{C}}{\mu} |\T^ix - \T^iy|\alf   \\
      & \leq &  \sum_{i=0}^{n-1} \frac{\hat{C}}{\mu} \frac{1}{\mu^{\alpha (n-i)}}
                                     |\T^nx - \T^ny|\alf
\end{eqnarray*}
\hfill $\Box$

\vspace{11 pt}

\subsubsection{Introduction of Holes}

A hole $H$ in [0,1] is a finite union of open intervals $H_l$, $1 \leq l \leq L$.
Let $I = [0,1]\backslash H$.  For each $n$, 
let $\displaystyle I^n = \cap_{i=0}^n T^{-i}I$ and
let $\displaystyle T^n = \T^n | I^n$. Let $K$ be the number of 
intervals of monotonicity for $T$ and let $\Q = \{ I_1, ... I_K \}$
be the partition of $I$ into those intervals.  
Denote by $I_{j,n}$ the part of $I_j$ which remains outside the hole for the 
first $n$ iterates of the map $\T$, i.e. $\displaystyle I_{j,n} = I_j \cap I^n$.
Note that $(T^n)'$ satisfies the same distortion bounds as $(\T^n)'$ wherever 
it is defined on each $\hat{I}^n_j$.

Let $d$ be the minimum length of the $I_j$ and $D$ their maximum length.
Define $C:= \C (2\delta)\alf$ and choose $\delta$ small enough that
\[
  \frac{2\alf(1+C)}{\mu\alf} < 1.
\]
We do this so that the tower we construct will satisfy the nonlinearity
condition (H2).  If necessary, reduce $\delta$ further so that $\delta \leq d$.
(If $\T$ is $C^2$, then we simply take $\delta = d$ since (H2) is not needed.)
$\delta$ represents the length scale of our reference intervals which we use to construct the tower $\Delta$.
Let $a := \max \left\{ \sqrt{\frac{2}{\mu}}, \frac{2\alf(1+C)}{\mu\alf} \right\}$.
Our sole condition on the hole is that its total measure must satisfy:

\vspace{11 pt}
\noindent
(A1) \hfill $\displaystyle mH \leq \frac{(1-a)^2 \mu \delta^2 
           \left( 1 - \sqrt{\frac{2}{\mu}}  \right) }{1+C}$.
      \hfill $\mbox{}$
\vspace{11 pt}

\noindent
This assumption ensures that the hole in the tower we construct satisfies condition (H$3'$).
(We do not need to satisfy (H3) since there will be no holes in the base of the tower we shall
construct.)  

Assumption (A1) also implies the following bound, which we will use to construct the tower,
\begin{equation}
h \leq \delta \frac{\mu - 2}{2}    \label{bound on hole}
\end{equation}
where $h$ is the maximum length of the intervals $H_l$.  To see that equation~(\ref{bound on hole})
is implied by (A1), we note that $\delta \leq \frac{1}{\mu}$ and write
\[
h \leq mH \leq \frac{(1-a)^2 \mu \delta^2 
           \left( 1 - \sqrt{\frac{2}{\mu}}  \right) }{1+C} \leq \delta \left( 1 - \frac{2}{\mu} \right) 
        < \delta \left( \frac{\mu}{2} -1 \right).
\]

\subsubsection{Statement of Results}
\label{expanding results}

We first state our result on the existence of an a.c.c.i.m.

\begin{theorem}
\label{expanding existence}
Let $\T$ be an expanding map of the interval with properties (a) and (b)
of Section~\ref{properties}.  Suppose $H$ is a hole satisfying (A1).
Then there exists a nontrivial absolutely continuous conditionally invariant measure for $T$
whose density is bounded above.  If $\alpha = 1$, the density is of bounded variation.

Suppose in addition that $T$ satisfies the following transitivity property:
for each interval $I_j$, there exists an integer $n_j$ such that
\[
I_j \cup T I_{j,1} \cup T^2 I_{j,2} \cup T^3 I_{j,3} \cup 
   \ldots \cup T^{n_j} I_{j,n_j} = I;  
\]
then the conditionally invariant density is bounded away from zero on I.
\end{theorem}

\noindent
{\bf Remark.}  We shall show in Section~\ref{accim bounds} that the eigenvalue of the a.c.c.i.m.\
satisfies the Lipschitz bound $1-\lambda \leq C_0 mH$ where $C_0$ depends only on the map $\T$
and the lengthscale $\delta$.  \\
\\
This theorem implies that Lebesgue measure has a well-defined exponential escape rate given 
by the a.c.c.i.m.  The escape rate relates to the question of uniqueness which is subtle
for open systems.  It is noted in \cite{chernov exp} that $C^2$ expanding maps may have families
of conditionally invariant densities of bounded variation which share the same support, but have
different eigenvalues, even when those maps are transitive on components and the holes are elements
of a finite Markov partition.  Nevertheless, since the spectrum of the transfer operator on $BV$ is stable
in the presence of small holes (\cite{liverani maume}), we are assured when $\alpha = 1$ that we
have found the unique density in $BV$ whose eigenvalue gives the escape rate of the reference 
measure.  When generalizing to $\alpha < 1$, our densities may no longer be in $BV$ and it is not
clear in what sense we can expect uniqueness, although we still have the correct escape rate
with respect to the reference measure.

Our next result concerns the convergence as $mH \rightarrow 0$  of the conditionally
invariant measures to an absolutely continuous invariant measure for $\T$ with no holes.

Let $\T$ and $H$ satisfy the hypotheses of Theorem~\ref{expanding existence}.
Since we are interested in arbitrarily small holes, we may further assume that 
either none of the closures of the holes $H_l$ contains an endpoint of one of the
$\hat{I}_j$, or if one does, that $H_l$ shrinks to this point in the
limit.  (If $H$ does not satisfy this property, we may take a smaller hole which does.)
In this way, we assure that no new intervals $I_j$ are created as we take the
limit $mH \rightarrow 0$.  In fact, each interval $I_j$ will either grow or remain
the same as the hole shrinks. 

We need the following assumption on the dynamics of $\T$.  We assume that for each 
interval $I_j$, there exists an $m_j$ such that
\begin{equation}
\label{eq:mixing}
I_j \cup \T I_j \cup \T^2 I_{j,1} \cup \T^3 I_{j,2} \cup 
   \ldots \cup \T^{m_j} I_{j,m_j-1} = [0,1].
\end{equation}
This means that the images of $I_j$ eventually cover [0,1], but we do not allow
a part that has fallen into the hole to be considered in the future - it is only
considered to cover that part of $H$ on which it lands.  This assumption is
obviously stronger than the transitivity assumption of Theorem~\ref{expanding existence}.

Note that this assumption on the dynamics of $\T$, in addition to the fact that 
$\T$ is expanding, implies that the map $\T$ with no holes has the covering property
and so has a unique a.c.i.m.\ whose density is bounded away from zero on [0,1].

\begin{theorem}
\label{convergence}
Let $\T$ and $H$ satisfy the hypotheses of Theorem~\ref{expanding existence}.  
Define $H_t = H$ and let $\{H_s \}$ for $s \in [0,t]$ be a sequence of holes 
with the following properties:
\begin{enumerate}
  \item[(1)] $mH_s \leq s$, $H_s \subset H_t$ and each component of $H_t$ contains at
    most one component of $H_s$;
  \item[(2)] either the closure of $H_t$ contains none of the endpoints of the
    $\hat{I}_j$, or if one endpoint is contained in the closure of $H_t$, then that endpoint is
    contained in the closure of $H_s$ for all $s \in [0,t]$;
  \item[(3)] equation (\ref{eq:mixing}) is satisfied by $\T$ acting on the intervals
    of monotonicity with respect to the hole $H_t$.
\end{enumerate}
Let $\nu_s$ denote the a.c.c.i.m.\ corresponding to $H_s$ obtained from Theorem~\ref{expanding existence}.
Then the $\nu_s$ converge weakly to $\nu$ as $s$ goes to zero, 
where $\nu$ is the unique absolutely continuous invariant measure for $\T$ with no holes.
If $\alpha = 1$, then the densities converge in $L^1([0,1])$ as well.
\end{theorem}

\subsection{Proof of Theorem \ref{expanding existence}}

Our proof takes the following steps:
\begin{enumerate}
   \item prove a growth lemma for intervals of length at least $\delta$;
   \item use the growth lemma to define a return time function and partition of $I$
     in order to build a tower;
   \item verify that the tower has the desired properties and satisfies
     conditions (H1)-(H$3'$) to conclude the existence of an a.c.c.i.m.;
   \item project the measure on the tower to a measure on $I$ and show it has
     the properties we claim.
\end{enumerate}

\subsubsection{A Growth Lemma}
\label{growth}

Although we have already fixed $\delta$ small enough to control the nonlinearities
in the tower, it should be noted that the results of Lemma~\ref{lem:stopping time}
hold for any $\delta'$ and $h$ satisfying
\[
  \frac{2h}{\mu -2} \leq \delta' \leq d
\] 
We will continue to use $\delta$ as the notation throughout this section, however.

Any interval $\Omega$ of length at least $\delta$ lying 
entirely in one of the $I_j$ will have at least one subinterval which
grows to cover one of the intervals of monotonicity exponentially fast.  
This is because
if $\T\Omega$ does not cover one of the $I_j$, then $\Q|\T\Omega$
can have at most two components, and $\T\Omega$ can intersect at
most one hole.  The measure of what does not fall into the hole is
\[
m(T(\Omega \cap I^1)) \geq \mu\delta - h \geq \mu\delta - \delta \left( \frac{\mu-2}{2} \right)
  = \delta \left( \frac{\mu + 2}{2} \right),
\]
where we have used equation~(\ref{bound on hole}) in the second inequality.
So at least one of the two components of $\Q|\T\Omega$ must have
length greater than or equal to $\delta(\frac{\mu+2}{4}) > \delta$
since $\frac{\mu +2}{4} > 1$.  In fact, if we follow this component,
call it $\omega_1$, we see that its image must also either cover one 
of the $I_j$ or else $\Q|\T\omega_1$ has at most two components, one 
of whose lengths is at least $\delta(\frac{\mu + 2}{4})^2$.  Always
following the larger component, we obtain a sequence of intervals 
$\omega_n$ of length
at least $\delta(\frac{\mu + 2}{4})^n$ which must eventually cover
one of the intervals of monotonicity of $T$ in exponential time depending
only on $\mu$, $\delta$ and $D$.

\begin{lemma}
\label{lem:stopping time}
Let $\Omega$ be any interval of length at least $\delta$ lying
entirely in one of the $I_j$.  There exists a countable partition $\Z$
of $\Omega$ and a stopping time $S(x)$ such that for each
$\omega \in \Z$,
\begin{enumerate}
  \item[(a)]  $S$ is constant on $\omega$;
  \item[(b)]  $\T$ is $C^{1+\alpha}$ on $\T^S \omega$;
  \item[(c)]  Either $\omega \subset I^S$ and $T^S \omega = I_j$ for some $1 \leq j \leq K$ or
    $\omega \subset I^{S-1}$ and $\T^S \omega \subset H$;
  \item[(d)] $m\{x \in \Omega \; | \; S(x)>n\} \leq D (\frac{2}{\mu})^n$,
    i.e. the elements of $\Z$ either fall into $H$ or grow at an exponential rate;
  \item[(e)]  $m\{ x \in \Omega \; | \; \T^{S(x)}(x) \in H \}
    \leq \frac{mH}{\mu - 2}$.
\end{enumerate}
\end{lemma}

\noindent
{\em Proof}.  Let $\Omega_0 := \Omega$.  Given $\Omega_{n-1}$,
which consists of a finite number of subintervals of $\Omega$,
we show how to form the set $\Omega_n$ inductively.  $\Omega_n$
represents those points in $\Omega$ which have not yet been
assigned a stopping time by time $n$.

We proceed one element of $\Omega_{n-1}$ at a time.  Fix
$\omega \in \Omega_{n-1}$.  We consider two cases depending
on $\T^n\omega$.

\vspace{11 pt}
\noindent
Case 1:  $\T^n\omega$ does not cover any of the $I_j$.\\
Then $\Q|\T^n\omega$ has at most two components and
$\T^n\omega$ intersects at most one $H_l$. Set $S(x) = n$ on
$\omega \cap \T^{-n}H$ and enter this interval as an element of the
partition $\Z$.  Note that this interval has length less than or
equal to $\frac{h}{\mu^n}$.  Put the (at most two) components of
$\T^{-n}\Q|\omega$ into the new set $\Omega_n$.  The stopping time
$S$ has not yet been defined on the elements of $\Omega_n$.

\vspace{11 pt}
\noindent
Case 2:  $\T^n\omega$ covers at least one of the $I_j$.\\
Set $S(x)=n$ on each of the components of $\T^{-n}\Q|\omega$ whose image
equals one of the $I_j$.  Each of these components
is entered as an element of the partition $\Z$.  $\T^n\omega \cap H$
has at most $L$ components.  Set $S(x)=n$ on each of the intervals
$\omega \cap \T^{-n}H_l$ and enter these as elements of the partition
$\Z$.  Note that $m(\omega \cap \T^{-n}H) \leq \frac{mH}{\mu^n}$.
There are at most two intervals in $\T^{-n}\Q|\omega$ left over where $S$ has
not been defined.  Include these intervals as elements of $\Omega_n$. 
Note that their images under $T^n$ each have length less than $D$.

\vspace{11 pt}
By construction, $\Omega_n$ has at most $2^n$ components, each
of length less than $ \frac{D}{\mu^n}$.  This implies
\[
  m \{ x\in \Omega \; | \; S(x)>n \} = m \Omega_n 
  \leq D \left( \frac{2}{\mu} \right)^n
\]
so that (d) is satisfied.  (a),(b) and (c) are clear from the construction
of $\Omega_n$.  

To check (e), we note that the most Lebesgue measure that can fall
into $H$ at the $n^{th}$ step is $\frac{mH}{\mu^n}$ from each
component of $\Omega_{n-1}$.  This means at most $mH \frac{2^{n-1}}{\mu^n}$
falls in at each step.  So the total measure of those elements of $\Z$
which fall in the hole before they can grow to cover one of the intervals of monotonicity is
\[
 m\{x\in \Omega \; | \; \T^{S(x)}(x) \in H \}  \leq
 \sum_{n=1}^{\infty} \frac{mH}{\mu} \left( \frac{2}{\mu} \right)^{n-1}
 = \frac{mH}{\mu - 2}
\]
\hfill $\Box$

\subsubsection{Building the Tower}
\label{construction}

We show how to build a tower with multiple bases from the growth
lemma.  We partition $I$ into $N$ intervals of length between
$\delta$ and $2\delta$.  Call them 
$\Lambda^{(1)}, \ldots \Lambda^{(N)}$.  Since each $\Li$ has length at 
least $\delta$, each will have a partition $\Z^{(i)}$ and a stopping 
time $S^{(i)}$ with properties (a)-(e) of Lemma~\ref{lem:stopping time}.

Proceeding one $\Li$ at a time, we define below the return time function $R$ 
and the partition $\Z$ of $I$ on whose elements $R$ is constant.
This will define the tower $\hat{\Delta}$.  Fix
$k$ between 1 and $N$.  On those $\omega \in \Z^{(k)}$ such that
$\T^{S^{(k)}(\omega)} \omega \subset H$, we set $R=S^{(k)}$ and include 
$\omega$ as an element of $\Z$.

For those $\omega \in \Z^{(k)}$ with $\T^{S^{(k)}(\omega)} \omega = I_j$ 
for some $j$, we set $R=S^{(k)}$, but do not include $\omega$
as an element of $\Z$ directly.  $I_j$ consists of a number of the 
$\Li$.  We partition $\omega$ into elements of the form 
$\omega \cap T^{-S^{(k)}(\omega)}\Li$ and include each of these as
elements of the partition $\Z$.  We do this because each of the
$\Li$ may have varying lengths and so their projections from the tower
will have derivatives which vary accordingly.  We subdivide the partition
$\Q$ into the intervals $\Li$ in order to have control over distortion bounds on the tower.

We arrive at a return time function $R$ defined on $I$ and a countable 
partition $\Z$ which respects the boundaries of   
$\Lambda^{(1)}, \ldots \Lambda^{(N)}$ and satisfies the following properties for
each $\omega \in \Z$,
\begin{enumerate}
  \item[(a)]  $R$ is constant on $\omega$;
  \item[(b)]  $\T$ is $C^{1 + \alpha}$ on $\T^R \omega$;
  \item[(c)]  Either $\T^R \omega \subset H$ or $\T^R \omega = \Li$
            and $T^R$ is defined on $\omega$;
  \item[(d)]  $m \{ x \in I \; | \; R(x)>n \} \leq N D (\frac{2}{\mu})^n$, 
            i.e. the elements of $\Z$ either fall into $H$ or return at 
            an exponential rate;
  \item[(e)] $m \{ x \in I \; | \; \T^{R(x)}x \in H \} 
           \leq N \frac{mH}{\mu -2}$
\end{enumerate}
Lemma~\ref{lem:distortion} implies that on each $\omega \in \Z$, we have the following 
distortion bound for $x,y \in \omega$ and $n \leq R(\omega)$.
\begin{equation}
\label{eq:Tdistortion} 
\left| \frac{(\T^n)'(x)}{(\T^n)'(y)} - 1 \right| \leq \C |\T^nx - \T^ny|\alf
\end{equation}

We identify our reference intervals $\Lambda^{(1)}, \ldots, \Lambda^{(N)}$ with $N$ disjoint intervals of unit length, 
$\hat{\Delta}_0^{(1)}, \ldots, \hat{\Delta}_0^{(N)}$, as follows.  Let 
$\hat{\Delta}_0 = \bigcup_{i=1}^{N} \hat{\Delta}_0^{(i)}$ be as in Section~\ref{multiple bases}.
Define a bijection $\pi_0 : \hat{\Delta}_0 \rightarrow I$ as piecewise linear with
$\pi_0( \hat{\Delta}_0^{(i)}) = \Li$ and
$\pi_0' | \hat{\Delta}_0^{(i)} = | \Li |$.  This makes $\delta \leq \pi_0'
\leq 2 \delta$ on $\hat{\Delta}_0$.  The partition $\Z$ and
return time function $R$ on $I$ induce a partition and return time function on
$\hat{\Delta}_0$.  We call these $\Z$ and $R$
as well: they are the partition and return time function introduced in Section~\ref{multiple bases}.
We define as before
\[ 
  \hat{\Delta} = \{ (x,n) \in \hat{\Delta}_0 \times \mathbb{N} : R(x) > n \}.
\]
Recall the notation $\hat{\Delta}_l = \hat{\Delta}|_{n=l}$ for
the $l^{th}$ level of the tower, and let $\hat{\Delta}_l^{(i)}$
denote the $l^{th}$ level above $\hat{\Delta}_0^{(i)}$.
If $\T^R \omega \subset H$, then we put a hole 
$\h_{R(\omega),j}$ in the $R(\omega)$ level of the tower
above $\pi_0^{-1}\omega$ in $\hat{\Delta}_0$.  This defines
the hole $\h$ in the tower $\hat{\Delta}$.  
Let $\Delta = \hat{\Delta} \backslash \h$ and in general 
$\Delta_l^{(i)} = \hat{\Delta}_l^{(i)} \backslash \h$.  Note that there are no holes in the base.

\subsubsection{Properties of the Tower}
\label{tower properties}

Let $\hat{F}$ be the tower map, $\F : \Delta \rightarrow \hat{\Delta}$.  Define a projection 
$\pi : \hat{\Delta} \rightarrow [0,1]$ such that 
$\pi \circ \F = \T \circ \pi$ on $\Delta$.  The elements of the Markov 
partition $\dilj$ are the maximal intervals on $\hat{\Delta}_l^{(i)}$
which project onto the dynamically defined
elements of $\Omega_l$ in the construction of the return time function
$R$ above the reference interval $\Li$.

For $x \in \Delta_0$, we have 
$\pi' \circ \F^l(x) \cdot (\F^l)'(x) = (\T^l)'\circ \pi(x) \cdot \pi'(x)$
and if $R(x)>l$ then $(\F^l)'(x) = 1$ so that
\begin{equation}
  \pi' \circ \F^l(x) = (\T^l)'\circ \pi(x) \cdot \pi'(x).
  \label{eq:projection}
\end{equation}

Now let $z \in \Delta_0^{(i)}$ with $R(z) = l + 1$.  Let $\F^lz = x$.
Then $\pi' \circ \F^{l+1}(z) \cdot (\F^{l+1})'(z) = 
(\T^{l+1})'\circ \pi(z) \cdot \pi'(z)$.  But $\F^{l+1}(z) \in \Delta_0^{(j)}$
so that $\pi' \circ \F^{l+1}(z) = |\Lambda^{(j)}|$ and 
$(\F^{l+1})'(z) = \F'(\F^l z) \cdot (\F^l)'(z) = \F'(x)$.  This yields
\begin{equation}
\F'(x) = (T^{l+1})'(\pi z) \frac{|\Li|}{|\Lambda^{(j)}|} \geq \frac{\mu^{l+1}}{2}.
                                                      \label{eq:derivative}
\end{equation}

We derive a distortion estimate for $\F$.  Let $x,y \in \dklj$ be such that $\F(x), \F(y) \in \hat{\Delta}_0^{(i)}$.
Using equations~(\ref{eq:Tdistortion}) and (\ref{eq:derivative}) we have
\begin{eqnarray*}
  \left| \frac{\F'(x)}{\F'(y)} - 1 \right|  & = & 
       \left|  \frac{(\T^{l+1})'(\pi \circ \F^{-l}x)}
                    {(\T^{l+1})'(\pi \circ \F^{-l}y)} - 1 \right|   \\ 
    & \leq & \C |\T^{l+1}(\pi \circ \F^{-l}x) 
          -  \T^{l+1}(\pi \circ \F^{-l}y)|\alf     \\
    & = & \C | \pi \circ \F (x) - \pi \circ \F (y) |\alf  \\
    & \leq & \C (2 \delta)\alf |\F (x) - \F (y)|\alf
\end{eqnarray*}

Let $F= \F|(\Delta \backslash \F^{-1}\h)$.  $F$ also satisfies the relation
$\pi \circ F = T \circ \pi$ on its domain, and so the above estimates hold for
$F'$.

\subsubsection{$\mathbf{(\Delta, F, m)}$ Satisfies Conditions (H1)-(H3$\mathbf{'}$)}

Recall properties (a)-(c) required of the tower map in Section~\ref{regularity of F} as well
as assumptions (H1)-(H3$'$) of Section~\ref{tower section}.
It is clear from the discussion of Section~\ref{tower properties} that $\F$ has properties (a)-(c) 
of the tower map with
$\gamma \geq \frac{\mu}{2}$, $\beta = \log \mu$, $C = \C(2\delta)\alf$,
and the same H\"{o}lder exponent $\alpha$ as the map $\T$.

(H1) is satisfied with $\theta = \frac{2}{\mu}$.  This is because
$m \hat{\Delta}_l = m \Delta_l + m \h_l$ and we have bounds on these
two quantities from the proof of Lemma~\ref{lem:stopping time}.
We have $ m \Delta_l \leq \frac{ND}{\delta} (\frac{2}{\mu})^l$ and
$  m \h_l \leq \frac{NmH}{\delta \mu}(\frac{2}{\mu})^{l-1}$.
This yields the desired estimate $m \hat{\Delta}_l \leq \frac{N}{\delta} (\frac{2}{\mu})^l$.

(H2) is satisfied by choosing $\delta$ small enough.

Choose $\xi = \min \{ \frac{1}{2} \log (\frac{\mu}{2}), \alpha \log \mu \}$.  Then
$e^{\xi} \leq \sqrt{\frac{\mu}{2}}$ so (H$3'$) becomes
\[
\sum_{l=1}^{\infty} \left( \frac{\mu}{2} \right)^{\frac{l-1}{2}}m\h_l
     < \frac{(1-a)^2}{1+C}.
\]
But $m\h_l \leq \frac{NmH}{\delta \mu}(\frac{2}{\mu})^{l-1}$.
So
\begin{eqnarray*}
\sum_{l=1}^{\infty} \left( \frac{\mu}{2} \right)^{\frac{l-1}{2}}m\h_l 
     & \leq & \sum_{l=1}^{\infty} \frac{NmH}{\delta \mu}
                   \left( \frac{2}{\mu}\right)^{l-1} \left( \frac{\mu}{2}\right)^{\frac{l-1}{2}} \\
     & = & \frac{NmH}{\delta \mu} \frac{1}{1 - \sqrt{\frac{2}{\mu}}}
\end{eqnarray*}
Using this estimate, we see that (H$3'$) will be satisfied if $H$ satisfies
\[
   mH < \frac{(1-a)^2}{1+C} \cdot
             \frac{\mu \delta (1 - \sqrt{\frac{2}{\mu}})}{N}
\]
which is slightly weaker than our assumption (A1).

\subsubsection{Existence and properties of an a.c.c.i.m.}
\label{accim bounds}

Since $(\Delta, F, m)$ satisfies (H1)-(H$3'$), we conclude using Theorem~\ref{tower theorem} that
there exists a $\varphi \in X$ such that $d \tilde{\nu} := \varphi dm$ is a nontrivial
a.c.c.i.m.\ with respect to $F$ acting on $\Delta$.  Let $\lambda$ be the
eigenvalue of $\tilde{\nu}$.  Now define
a measure $\nu$ on $I$ by 
\[
  \nu(A) := \tilde{\nu}(\pi^{-1}A)
\]
for any Borel subset $A$ of $I$.  Then it is clear that $\nu$ will be
conditionally invariant with respect to $T$ with the same eigenvalue
$\lambda$ since for any Borel $A \subset I$,
\[
\begin{array}{rcccl}
 \nu(T^{-1}A) & := & \tilde{\nu}(\pi^{-1} \circ T^{-1}A) & = &  \tilde{\nu}(F^{-1} \circ \pi^{-1}A) \\
                   & = & \lambda \tilde{\nu}(\pi^{-1}A)  & =: &  \lambda \nu (A).
\end{array}
\]

We show that $\nu$ is absolutely continuous with respect to Lebesgue measure 
with a density that is bounded above.  For this we will need to estimate the number
of preimages that a Borel set $A$ can have on each level of the tower.  We consider
the tower $\Delta^{(i)}$ above one of the $\Delta_0^{(i)}$.  From the proof of 
Lemma~\ref{lem:stopping time} we know that the number of partition elements 
$\Delta_{l-1,j}^{(i)}$ which move up to $\Delta_l^{(i)}$ is at most $2^{l-1}$ (these
are the preimages under $\pi$ of the elements of $\Omega_{l-1}$).  Although each
of these elements may be split into 2 or more pieces on $\Delta_l^{(i)}$, 
the map $T$ is monotonic on each $\pi(\Delta_{l-1,j}^{(i)})$ and so 
$T \circ \pi(\Delta_{l-1,j}^{(i)})$ covers $A$ at most once.  Thus 
$\pi(\Delta_l^{(i)})$ covers $A$ at most $2^{l-1}$ times.

This observation and equation~(\ref{eq:projection}) yield the following estimate.
\begin{eqnarray*}
 \nu(A) & =    & \sum_{k=1}^{N} \sum_{l=0}^{\infty} \sum_j
                   \int_{\pi^{-1}A \cap \dklj} \varphi dm   \\
        & \leq & \sum_{k=1}^{N} \sum_{l=0}^{\infty} \sum_j 
                   \sup_{\dklj} \varphi \; m(\pi^{-1}A \cap \dklj) \\
        & \leq & \sum_{k=1}^{N} \sum_{l=0}^{\infty} \sum_j
                   \frac{1}{\lambda^l} \sup_{\Delta_0^{(k)}} \varphi \; 
                   \frac{1}{\delta \mu^l} m(A \cap \pi \dklj)  \\
        & \leq & \sum_{k=1}^{N} \sum_{l=0}^{\infty}
                   \frac{1}{\lambda^l} \sup_{\Delta_0^{(k)}} \varphi \;
                   \frac{1}{\delta \mu^l} 2^{l-1} mA \\
        & \leq & \frac{N mA}{2 \delta} \sup_{\Delta_0} \varphi \sum_{l=0}^{\infty} 
                   \left( \frac{2}{\mu} \right)^{\frac{l}{2}}
\end{eqnarray*}
where we have used the fact that $\lambda \geq \sqrt{\frac{2}{\mu}}$ in the last line.
This proves that $\nu << m$ with bounded density.

Let $d\nu = \psi dm$. 
Choose a $\Delta_0^{(i)}$ such that $\inf_{\Delta_0^{(i)}} \varphi = \rho >0$. 
Then $\inf_{\Li} \psi \geq \frac{\rho}{2\delta}$.
If $T$ satisfies the transitivity property in the statement of
Theorem~\ref{expanding existence}, then there exists an $n_0$ such that the first $n_0$
images of $\Li$ under $T$ cover $I$.  Thus for any $x \in I$, there exists
an $n \leq n_0$ and $z \in \Li$ such that $T^nz = x$.
\[
  \lambda^n \psi(x) = \mathcal{P}^n \psi (x) := \sum_{T^ny = x} \frac{\psi(y)}{|(T^n)'(y)|}
                     \geq \frac{\psi(z)}{|(T^n)'(z)|}
\]
Let $\eta$ denote $\max_{[0,1]} |\T'|$.  Then 
\[
\inf_{I} \psi \geq \frac{\rho}{2 \delta \eta^{n_0}}.
\]
Thus $\nu$ is proportional to Lebesgue measure on $I$.

Since $\varphi \in X_M$, we have $|\mathcal{P}\varphi|_1 = \lambda \geq 1 - qM$ from the estimates
in Section~\ref{choice of M}.  Recall that $M = \frac{b}{1-a}$ where $a$ and $b$ depend only on the
map $\T$ and the length scale $\delta$.  
Similarly, $q = \sum_{\ell \geq 1} e^{\xi(\ell -1)} m\h_{\ell} \leq \frac{mH}{\delta^2(\mu - \sqrt{2\mu})}$.
From this we see that $\lambda \geq 1 - C_0 mH$, which is equivalent to the Lipschitz bound claimed
in the remark after Theorem~\ref{expanding existence}.
\\
\\
\noindent
{\bf Estimating the variation of $\mathbf{\varphi}$.}
If $\alpha =1$, we can use the bound on $\| \varphi \|_r$ to estimate the variation
of $\varphi$ since on $\Delta_l$, $|\varphi'| \leq \| \varphi \|_r \varphi e^{\xi l}$ 
by definition of $ \| \cdot \|_r$.  In fact, the following calculation shows that the projection
of any function in $X$ has bounded variation on $I$.

Let $\bigvee_J f$ denote the variation of $f$ on the interval $J$.  We estimate $\bigvee_I \psi$
using the following identity
\[
\psi(x) = \sum_{i,l,j} \frac{\varphi(\pilj x)}{\pi' (\pilj x)} 
\]
where $\pilj$ denotes the inverse of $\pi|_{\dilj}$.  We estimate
\[
\bigvee_I \psi \leq \sum_{i,l,j} \bigvee_{\pi(\dilj)} \varphi \circ \pilj \cdot (\pilj)'
      + 2 |\varphi \circ \pilj \cdot (\pilj)'|_{\infty}.
\]
The second term can be estimated using equation~(\ref{eq:projection}) and the fact that $\varphi \in X_M$,
\begin{equation}
\label{variation second term}
|\varphi \circ \pilj \cdot (\pilj)'|_{\infty} \leq Me^{\xi l} \cdot \frac{1}{\delta \mu^l} .
\end{equation}

To estimate the first term, we note that equation~(\ref{eq:projection}) and the distortion bounds of 
Lemma~\ref{lem:distortion} imply that
\[
\left| \frac{(\pilj)'(x)}{(\pilj)'(y)} - 1 \right| \leq \tilde{C} |x - y|
\]
which in turn implies $\left| \frac{(\pilj)''}{(\pilj)'} \right| \leq \tilde{C}$.
We use this to estimate
\begin{eqnarray}
\bigvee_{\pi(\dilj)} \varphi \circ \pilj \cdot (\pilj)' & \leq &
    \int_{\pi(\dilj)} |(\varphi \circ \pilj \cdot (\pilj)')'| \; dm  \nonumber \\
    & \leq &  \int_{\pi(\dilj)} | \varphi' \circ \pilj \cdot ((\pilj)')^2| \; dm  \nonumber  \\ 
    &      &  \; \; +  \int_{\pi(\dilj)} \varphi \circ \pilj |(\pilj)''| \; dm           \nonumber   \\
    & \leq & \frac{1}{\delta \mu^l} \int_{\dilj} |\varphi'| \; dm + \tilde{C} \int_{\dilj} \varphi \; dm 
                    \nonumber    \\
    & \leq & \frac{M e^{\xi l} }{\delta \mu^l} \int_{\dilj} \varphi \; dm + \tilde{C} \int_{\dilj} \varphi \; dm .
                    \label{variation first term}
\end{eqnarray}

Putting together equations~(\ref{variation second term}) and (\ref{variation first term}) we obtain
\begin{eqnarray}
\bigvee_I \psi & \leq & \sum_{i,l,j} \frac{M e^{\xi l} }{\delta \mu^l} \int_{\dilj} \varphi \; dm + 
   \tilde{C} \int_{\dilj} \varphi \; dm + 2 \frac{M e^{\xi l} }{\delta \mu^l} \nonumber \\
   & \leq & \tilde{C} + \sum_{i,l} \frac{3M e^{\xi l} }{\delta \mu^l } 2^{l-1}  \nonumber \\
   & \leq & \tilde{C} + \frac{3NM}{2 \delta} \sum_l \left( \frac{2}{\mu} \right)^{\frac{l}{2}} 
                        \label{uniform variation bound}
\end{eqnarray}
where we have used $e^{\xi} \leq \sqrt{\frac{\mu}{2}}$ in the last line.

\subsection{Proof of Theorem \ref{convergence}}

We prove under the conditions stated in Theorem~\ref{convergence}
that the conditionally invariant measures found
via the tower construction will converge to the unique absolutely continuous invariant
measure of the expanding map
with no holes as $mH \rightarrow 0$.  Our strategy will be to prove that the
sequence of conditionally invariant densities generated by the sequence of holes
is uniformly bounded above and below away from zero. The measures converge weakly to a 
density which satisfies the same bounds.  This is an invariant measure for the expanding map with
no holes and since it is equivalent to Lebesgue measure on [0,1], it must be unique due
to the transitivity property of $\T$ (equation~(\ref{eq:mixing})).  If $\alpha = 1$, we
prove $L^1$ convergence of the densities as well.

Let $\T$ and $\{ H_s \}$ be as in the statement of Theorem~\ref{convergence}.
Let $I_j^t$ be the
intervals of monotonicity created by $H_t$ and let $I^t = [0,1]\backslash H_t$.  

There exists a conditionally invariant measure $\varphi_t$ on the tower $\Delta$
constructed as in Section~\ref{construction} with eigenvalue $\lambda_t$ (we suppress
the dependence of $\Delta$ on $H_t$ in our notation).  This in 
turn induces a conditionally invariant density $\psi_t$ for 
$T_t := \T|I^t \cap \T^{-1}I^t$.  Let $\delta$ from
the proof be fixed and note that $N \leq \frac{1}{\delta}$ where $N$ represents
the number of $\Li$ needed to cover $I^t$.  The constant $\delta$ does not change
as the hole shrinks since the intervals of monotonicity of the map with holes can
only grow and no new ones are created.

We essentially need to repeat the arguments given in Section~\ref{accim bounds}, but
taking care to show that the upper and lower bounds are uniform in $s$.  This is
trivial in the case of the upper bound since $\delta$ is fixed and $N \leq \frac{1}{\delta}$.
Also $\sup_{\Delta_0} \varphi_s \leq M$ where $M= \frac{1+C}{1-a}$.  All the constants
here are independent of $s$.  This also holds for the estimate on $\bigvee_I \varphi_s$.

We need to do the lower bound more carefully.  We have
\[
   1  \; \;  =  \; \;   \sum_{l=0}^{\infty} \int_{\Delta_l} \varphi_t dm
      \; \; \leq \; \; \sup_{\Delta_0} \varphi_t \sum_{l=0}^{\infty} \frac{1}{\lambda_t^l}m \Delta_l
      \; \; \leq \; \; \sup_{\Delta_0} \varphi_t \sum_{l=0}^{\infty} \frac{1}{\lambda_t^l}
                       \frac{ND}{\delta} \left( \frac{2}{\mu} \right)^l  . 
\]
Since $\lambda_t \geq \sqrt{\frac{2}{\mu}}$ and $D<1$, this yields the inequality:
\[
\sup_{\Delta_0} \varphi_t \geq \delta^2 \left( 1 - \sqrt{\frac{2}{\mu}} \right)  =: \kappa
\]
So there exists a $\Delta^{(i)}_0$ such that the maximum of $\varphi_t$ on that
interval is at least $\kappa$.  The regularity of $\varphi_t$ then forces
$\inf_{\Delta^{(i)}_0} \varphi_t \geq \frac{\kappa}{1+M}$.
This in turn implies that
\[
 \inf_{\Li} \psi_t \geq \frac{\kappa}{1+M}. 
\]
Now since $\Li$ grows (exponentially fast) to cover one of the intervals $I_j^t$, and each
of these grows in at most $m_j$ iterates to cover all of [0,1], then there exists an
$m_0$ such that $\Li$ grows to cover all of [0,1] in at most $m_0$ iterates.  So if
$\mathcal{P}_t$ is the Perron-Frobinius operator associated with $T_t$, then
for any $x \in I^t$, there exists an $n \leq m_0$ and $z \in \Li$ such that $T_t^nz = x$.
\[
  \lambda_t^n \psi_t(x) = \mathcal{P}_t^n \psi_t(x) := \sum_{T_t^ny = x} \frac{\psi_t(y)}{|(T_t^n)'(y)|}
                     \geq \frac{\psi_t(z)}{|(T_t^n)'(z)|}
\]
Then 
\begin{equation}
\inf_{I^t} \psi_t \geq \frac{\kappa}{(1+M) \eta^{m_0}}.
\label{eq:lower bound}
\end{equation}

We are now ready to shrink the hole.  Note that once the constants
$\kappa$ and $M$ have been chosen for $t$, then they will work for each $s$ as well.
Also, $m_0$ can be chosen large enough to work for all $s \in [0,t]$ since the intervals
$I_j^s$ can only be larger than the intervals $I_j^t$.  So we have the same lower
bound from (\ref{eq:lower bound}) for all the conditionally invariant densities $\psi_s$:
\[
\inf_{I^s} \psi_s \geq \frac{\kappa}{(1+M) \eta^{m_0}}.
\]

Now let $\nu_s$ be the measures conditionally invariant with respect to $T_s$ having
densities $\psi_s$ and eigenvalues $\lambda_s$.  Note also that by the proof of the
existence of the conditionally invariant measures on the tower, we have 
$\lambda_s \rightarrow 1$ as $s \rightarrow 0$.
Now take any sequence $\{s_i\}$ with $s_i \rightarrow 0$.  The sequence of
measures $\{ \nu_{s_i} \}$ is precompact and so we may chose a subsequence 
$\{ \nu_{s_k} \}$ such that the $\nu_{s_k}$ converge weakly to a probability measure $\nu$ on
[0,1].  Now for any Borel subset $A$ of [0,1], we have
\[
  \frac{\kappa}{(1+M) \eta^{m_0}} m(A) \leq \nu(A) \leq \frac{M}{\kappa } m(A) 
\]
so that $\nu$ is equivalent to Lebesgue measure on [0,1].  Also
\[
\begin{array}{rcccl}
\nu(\T^{-1}A) & = & \lim_{k \rightarrow \infty} \nu_{s_k}(\T^{-1}A) 
          & = & \lim_{k \rightarrow \infty} \nu_{s_k}(T_{s_k}^{-1}(A \backslash H_{s_k}))  \\
              & = & \lim_{k \rightarrow \infty} \lambda_{s_k} \nu_{s_k}(A \backslash H_{s_k})
          & = &  \lim_{k \rightarrow \infty} \lambda_{s_k} \nu_{s_k}(A) \\
              & = & \nu(A)
\end{array}
\]
so that $\nu$ is an invariant measure for $\T$.
These facts, combined with the mixing property of $\T$, makes $\nu$ the unique absolutely
continuous invariant measure for $\T$.  Since $\nu$ is unique, the above argument forces 
the convergence of all the measures $\nu_s$ to $\nu$ as $mH_s \rightarrow 0$.

If $\alpha =1$, we use a similar argument on the densities $\psi_s$ since we have uniform bounds on their 
variations given by equation~(\ref{uniform variation bound}).  We extend $\psi_s$ to $L^1(\I)$ by simply setting $\psi_s = 0$ on $H_s$.  Since we have
uniform upper bounds on $\psi_s$ and the number of components of the hole does not increase as 
$s \rightarrow 0$, the variation of $\psi_s$ will still be finite and uniformly bounded in $s$.  Call this
bound on the variations $B$.

Since the $\psi_s$ lie in a compact subset of $L^1(\I)$, we can choose a subsequence $\{ \psi_{s_k} \}$
which converges in $L^1$ to a function $\psi$ whose variation is at most $B$.  $\psi$ also has the
same lower bound as the $\psi_s$.  It only remains to show that $\psi$ is an invariant density for
$\T$.

Let $I_s = \I \backslash H_s$ and let $I_s^1 = I_s \cap T_s^{-1} I_s$.
We have the relation $\mathcal{P}_s \psi_s = \lambda_s \psi_s$ for each $s$ where $\mathcal{P}_s$
is the Perron-Frobinius operator corresponding to $T_s$.  Let $\hat{\mathcal{P}}$ be the
Perron-Frobinius operator corresponding to $\T$ and note that
\[
\mathcal{P}_s f = \hat{\mathcal{P}} (f \cdot 1_{I_s^1}).
\]
Noting that $|\psi_s|_{\infty} \leq 2 + B$, we are ready to estimate,
\begin{eqnarray*}
\int_{\I} |\hat{\mathcal{P}} \psi - \psi | & \leq &
  \int_{\I} |\hat{\mathcal{P}} \psi - \hat{\mathcal{P}} \psi_s| + |\hat{\mathcal{P}} \psi_s - \psi_s|
  + |\psi_s - \psi|  \; dm \\
  & \leq & 2|\psi_s - \psi|_1 + | \hat{\mathcal{P}} \psi_s - \hat{\mathcal{P}} (\psi_s \cdot 1_{I_s^1})|_1 
  + | \hat{\mathcal{P}} (\psi_s \cdot 1_{I_s^1}) - \psi_s |_1 \\
  & \leq & 2|\psi_s - \psi|_1 + |\psi_s \cdot 1_{\I \backslash I_s^1}|_1 + |\mathcal{P}_s \psi_s - \psi_s|_1  \\
  & \leq & 2|\psi_s - \psi|_1 + (2+B) m (\I \backslash I_s^1) + (1 - \lambda_s ).
\end{eqnarray*}
This last line approaches zero for $s = s_k$ as $k \rightarrow \infty$.  We conclude that
$\hat{\mathcal{P}} \psi = \psi$ so that $\psi$ is the density for the unique SRB measure
for $\T$.  Since $\psi$ is unique, the above argument forces the convergence in $L^1$ of all the densities
$\psi_s$ to $\psi$ as $mH_s \rightarrow 0$.

\vspace{11 pt}
\noindent
{\bf Acknowledgments.}  The author is grateful to L.-S. Young for many helpful discussions and guidance during
the writing of this paper.

\end{document}